\documentclass[11pt]{article}
\usepackage{amssymb, amsmath, amsthm}

\begin{document}
\title{Laplacian Solitons and Symmetry in $G_2$-geometry}
\author{\textsc{Christopher Lin}}
\date{}

\maketitle

\hskip 5cm \footnotesize\textsc{Abstract}\\
 In this paper, it is shown that (with no additional assumptions) on a compact $7$-dimensional manifold which admits a $G_2$-structure soliton solutions to the Laplacian flow of R. Bryant 
 can only be shrinking or steady.  We also show that the space of symmetries (vector fields that annihilate via 
 the Lie derivative) of a torsion-free $G_2$-structure on a compact $7$-manifold is canonically isomorphic to $H^1(M,\mathbb{R})$.  Some comparisons with Ricci solitons are also discussed, along with some future directions of exploration.

\normalsize
\section{Introduction}
 Let $M$ be a $7$-dimensional manifold
that admits a $G_2$-structure given by a non-degenerate $3$-form $\varphi$.  A natural geometric flow when $M$ is 
compact is the
Laplacian flow first suggested by R. Bryant in \cite{RB}:
\begin{equation}\label{lap}
 \frac{\partial \varphi}{\partial t} = -\Delta_{\varphi}\varphi
\end{equation}
for a family $\varphi = \varphi(t)$ of $G_2$ structures, where
$\Delta_{\varphi}$ denotes the Hodge Laplacian with respect to the
metric induced by $\varphi(t)$ (hence nonnegative-definite).\footnote{The sign in front of the Laplacian 
turns out to be purely a technical convention, as explained in \cite{KMT}.} 
 The original intention of the equation
(\ref{lap}) is to flow $\varphi$ to a torsion-free $G_2$-structure,
 since $\varphi$ being torsion-free is equivalent to being harmonic with respect to the metric
$g_\varphi$ it
induces.  Although it is not clear what it means in this case, when $M$ is not compact the flow (\ref{lap}) 
still makes sense.  In fact, when $M$ is not compact we suspect that a more general flow is needed (see \cite{SK} 
for some general results in this direction).  The short-time existence and uniqueness of (\ref{lap}) for a closed 
initial $G_2$-structure when $M$ is compact have been established in \cite{BX} and \cite{XY}.  \\

As in the Ricci flow, let us consider solutions of the form
\begin{equation}\label{soliton}
\varphi(t) = \tau(t)f_t^*\varphi
\end{equation}
for a fixed $G_2$-structure $\varphi$, pulled back by a
smooth family family $f_t$ of diffeomorphisms of $M$, and where
$|\tau(t)|>0$ is a scaling factor. Such special solutions are
called \text{\it{solitons}}, and in the present context of $G_2$ geometry it has 
already appeared in \cite{KMT} by Karigiannis, McKay, and Tsui.  We will call solutions of the
form (\ref{soliton}) Laplacian solitons, or just solitons if no
confusion arises.  Just as in the Ricci flow, we seek a static description of a soliton by 
substituting (\ref{soliton}) into (\ref{lap}), and we will arrive at the equation 
\begin{equation}\label{sol}
 \rho\varphi \,+\, L_{X}\varphi = -\Delta\varphi
\end{equation}  
for some constant $\rho$ and vector field $X$.  Thus equivalently, we can define a (Laplacian) soliton to be a $G_2$-structure $\varphi$ that satisfies (\ref{sol}).  \\  

Analogous to Ricci solitons, we can define:

\newtheorem{definition}{Definition}
\begin{definition}
 Let $\varphi$ be a $G_2$-structure, and $X$ a vector field
 on $M$.  We say that $\big(\varphi,X\big)$ is a Laplacian soliton
 if equation (\ref{sol}) is satisfied for some constant $\rho\in\mathbb{R}$.
   We say $\big(\varphi,X\big)$ is an expanding soliton if $\rho
   >0$, a steady soliton if $\rho = 0$, and a shrinking soliton if
   $\rho <0$.
\end{definition}
We would like to point out that a torsion-free $G_2$-structure is steady, and is the most trivial
example of a soliton.  Also, just to distinguish between the dynamic and static versions of the soliton concept, 
we will refer to (\ref{soliton}) as a \text{\it{soliton solution}} whereas the terminology of \text{\it{soliton}} 
will be reserved for the definition above.\\

Work on a dual equation (\text{\it{the Laplcian coflow}}) to (\ref{lap}) and the corresponding soliton equation analogous to (\ref{sol}) have been 
done in \cite{KMT}.  The original intention in \cite{KMT} was to study (\ref{sol}), but instead they focused 
on the coflow version because there is available a special cohomogeneity-$1$ ansatz.  In this paper we focus on a 
more detailed examination of the fundamental equation (\ref{sol}).  In addition, we found out very recently 
in \cite{WW2} that Weiss and Witt had also studied soliton solutions to the $L^2$-gradient flow of a Dirichlet-type functional they proposed in an earlier paper \cite{WW}.  In \cite{WW2} similar results to the ones in this paper and in \cite{KMT} appeared, but note that the $L^2$-gradient flow of their energy functional is a different equation from the Laplacian flow (or coflow).     \\

One of the main results of this paper is Corollary \ref{main1}, which says that there are no compact expanding solitons,
 and no compact steady solitons except torsion-free $G_2$-structures.  A similar result was proved in \cite{KMT}, where 
 the soliton was assumed to be \text{\it{coclosed}} (and \text{\it{closed}} for the original flow (\ref{lap})).  Our result 
 does not depend on the soliton $\varphi$ being closed, and the proof follows directly from a fundamental identity established in Lemma \ref{wonder} in the same section.  Although the short-time existence of the Laplacian flow 
 (and of the coflow as well) have only been established for closed/coclosed structures, Corollary \ref{main1} is still 
 valuable because we think (\ref{sol}) is still an interesting equation in its own right.  \\
 
 The other main result of this paper is Corollary \ref{isom}, which says that for a torsion-free $G_2$-structure $\varphi$ 
 on a compact $7$-manifold, the space of all vector fields $X$ such that $L_X\varphi = 0$ is isomorphic to 
 $H^1(M,\mathbb{R})$.  We will refer to such a vector field as a 
 \text{\it{symmetry}} of the $G_2$-structure $\varphi$, for brevity.  We view 
 Corollary \ref{isom} as a kind of rigidity result, because for $\varphi$ fixed the soliton equation (\ref{sol})
 is invariant only by adding such vector fields.\footnote{We shall see that 
 this is again rooted in Lemma \ref{wonder}.}  Our original goal was to do the same for 
 compact shrinking solitons as well, but at the present we do not have such an analogous result.

\section{$G_2$-Structures and Torsion Forms}
\quad We give a brief review of the background relating to $G_2$-structures here.  The standard reference for this is 
the book \cite{DJ} by Dominic D. Joyce, although the papers \cite{RB} and \cite{SK} are also good sources. \\
 
The group $G_2$ is a compact, connected, simply-connected Lie group sitting in $SO(7)$.  Algebraically, it can be defined 
as the Automorphism group of the Octonians.  It appears as one of the exceptional holonomy groups in the classification by Berger, et al.  The working definition of $G_2$ that we will adopt is the following.  Consider the differential $3$-form 
\begin{equation}\label{natural}
 \varphi_0 = dx_{123} + dx_{145} + dx_{167} + dx_{246} - dx_{257} - dx_{347} - dx_{356}
\end{equation}
in $\mathbb{R}^7$, where $dx_{ijk} = dx_i\wedge dx_j \wedge dx_k$.  The group $G_2$ can be defined as the subgroup in $GL(7,\mathbb{R})$ that preserves $\varphi_0$, which 
means \[
  \varphi_0 \big(g(v_1), g(v_2), g(v_3)\big) = \varphi_0\big(v_1, v_2, v_3\big) 
\] 
for all vectors $v_1, v_2, v_3 \in \mathbb{R}^7$ and every $g\in G_2$.  The peculiar form that (\ref{natural}) takes 
reflects the combinatorial nature of permuting the Octonians.  \\   

From a principal bundle point of view, a $G_2$-structure on a $7$-dimensional manifold $M$ is just a sub-bundle with structure group $G_2$, of the $GL(7,\mathbb{R})$-frame bundle over $M$.  In other words, one can find local frames of 
the tangent bundle $TM$ such that all the transition functions value in the group $G_2$.  Because $G_2 \subset SO(7)$, a 
$G_2$-structure induces an orientation and a unique metric (which we will write $g_\varphi$) on $M$.  Take any local frame $\{e_i\}_{i=1}^7$ in a 
$G_2$-structure, then the $3$-form 
\begin{equation}\label{3form}
 \varphi = \omega_{123} + \omega_{145} + \omega_{167} + \omega_{246} - \omega_{257} - \omega_{347} - \omega_{356}
\end{equation} 
is well-defined over all of $M$, where $\omega_i$ is the local dual $1$-form to $e_i$.  
Conversely, if a $3$-form $\varphi$ on $M$ can be locally represented as (\ref{3form}) with respect to a frame, then the 
transition functions of such frames value in $G_2$ and we have a $G_2$-structure.  As a result, a $G_2$-structure is 
equivalent to a $3$-form $\varphi$ locally represented as in (\ref{3form}), and this is what we will refer to as a 
$G_2$-structure.  It is well-known that a $G_2$-structure exists if and only if the $7$-manifold is orientable and spin.  \\   

A $G_2$ structure $\varphi$ induces a point-wise orthogonal
decomposition (with respect to $g_\varphi$) of $p$-forms on $M$:
\begin{align}
 \Omega_7^2 &= \{X\lrcorner\varphi\,|\, X\in \Gamma(TM)\} = \{\beta\in\Omega^2\,|\,\ast(\varphi\wedge\beta)=-2\beta\}\notag\\
 \Omega_{14}^2 &=
 \{\beta\in\Omega^2\,|\,\beta\wedge\ast\varphi=0\}= \{\beta\in\Omega^2\,|\,\ast(\varphi\wedge\beta)=\beta\}\notag\\
 \Omega_1^3 &= \{f\varphi\,|\,f\in C^{\infty}(M)\}\notag\\
 \Omega_7^3 &=
 \{X\lrcorner\ast\varphi\,|\,X\in\Gamma(TM)\}\notag\\
 \Omega_{27}^3 &= \{h_{ij}g^{jl}dx_i\wedge
 \big(\frac{\partial}{\partial x^l}\lrcorner\varphi\big)\,|\,
 h\in Sym^2(T^*M), Tr_g(h)=0\},\notag
\end{align}
where
\begin{align}
 \Omega^2 &= \Omega_7^2 \oplus \Omega_{14}^2  \notag\\
 \Omega^3 &= \Omega_1^3 \oplus \Omega_7^3 \oplus \Omega_{27}^3.
 \notag
\end{align}
Then we can write
\begin{align}\label{tor}
 d\varphi &= \tau_0\ast\varphi + 3\tau_1\wedge\varphi +
 \ast\tau_3\notag\\
 d\ast\varphi &= 4\tau_1\wedge\ast\varphi + \ast\tau_2,\notag
\end{align}
where $\tau_0\in\Omega_1^0, \tau_1\in\Omega^1_7,
\tau_2\in\Omega^2_{14}$, and $\tau_3\in\Omega^3_{27}$ are called
the \text{\it{torsion forms}}.  The fact that the same $1$-form $\tau_1$ appears in
the decompositions of $d\varphi$ and $d\ast\varphi$ is non-trivial, but can
be shown via some computations (see \cite{SK}).  The terminology of torsion forms comes from the following.  
Having a $G_2$-structure $\varphi$ does not mean the holonomy group of $g_{\varphi}$ is contained in $G_2$.  The 
additional condition that is needed is the so-called \text{\it{torsion-free}} condition.  We say that a $G_2$-structure 
is torsion-free if $\varphi$ solves the nonlinear system of partial differential equations $\nabla\varphi = 0$, where 
$\nabla$ is the covariant derivative induced by $g_\varphi$.  It was shown in \cite{FG} that a $G_2$-structure is torsion-free 
if and only if it is closed and co-closed (with respect to the hodge star induced by $g_\varphi$).\footnote{This is an 
entirely local property, it is independent of whether or not $M$ is compact.}  Thus when $M$ is compact, $\varphi$ being 
torsion-free is equivalent to it being harmonic with respect to $g_\varphi$.  In view of the torsion forms defined 
above, we see that $\varphi$ is torsion-free if and only if  all four torsion forms
vanish on $M$. A $7$-manifold $M$ that admits a torsion-free $G_2$-structure has its Riemannian holonomy (with respect to 
$g_{\varphi}$) a subgroup of $G_2$, and such manifolds are simply known as \text{\it{$G_2$ manifolds}}.\\   

\section{The Soliton Equation}\label{rambling}
From direct computations, we see that a soliton
solution (\ref{soliton}) to the Laplacian flow (\ref{lap}) satisfies

\begin{equation*}
 \dot{\tau(t)}f_t^*\varphi \,+\, \tau(t)f_t^*\big(L_{X(t)}\varphi\big)  = -\tau(t)^{1/3}f_t^*\big(\Delta\varphi\big),
\end{equation*}
which is exactly (\ref{sol}) and where the vector field $X(t)$ is the infinitesimal generator of
the diffeomorphism $f_t$.  We have also used the following
fact:

\newtheorem{lemma}{Lemma}
\begin{lemma}
 If $\varphi$ is a $G_2$-structure, then any non-zero constant
 multiple $c\varphi$ is also a $G_2$-structure, and $g_{c\varphi}
 = c^{2/3}g_{\varphi}$.\footnote{$c>0$ preserves the orientation given by $\varphi$,
 $c<0$ reverses it.}
\end{lemma}

 Note that we dropped the subscript on the Laplace operator
since $\varphi$ is now fixed. Then we see that there is a soliton
solution (\ref{soliton}) to the flow (\ref{lap}) only if

\begin{equation*} 
\rho\varphi \,+\, L_{X}\varphi = -\Delta\varphi
\end{equation*}
for some vector field $X$ on $M$, where we have frozen at a time
$t$ and the constant $\rho = \dot{\tau}/\tau^{1/3}$.  As in the
case of Ricci Solitons, one can show that given a vector field $X$
and a $G_2$ structure $\varphi$ that satisfy (\ref{sol}) one can
generate a solution of the form (\ref{soliton}) to (\ref{lap}). \\

 We would also like to point out that the
soliton equation (\ref{sol}) is scale-invariant in the following sense.  Note
that given any $G_2$-structure $\varphi$ satisfying (\ref{sol}),
then for any $c\varphi$, $c\ne 0$, equation (\ref{sol}) is again
satisfied for $\tilde{\rho} = c^{-2/3}\rho$ and
$\tilde{X}=c^{-2/3}X$.

\section{Compact Solitons}\label{cpt}
 In this section we show that

\begin{lemma}\label{wonder}
 Let $M$ be a compact $7$-manifold.  For any $G_2$-structure $\varphi$ on
 $M$, vector field $X$, and $f\in C^{\infty}(M)$, we have
 \begin{equation*}
   \int_M L_X\varphi\wedge\ast f\varphi = -3\int_M df\wedge\ast
   X^{\flat}
 \end{equation*}
\end{lemma}

\text{\bf{Proof.}}   We have \[
    L_X\varphi = X\lrcorner d\varphi + d(X\lrcorner\varphi).
\]

From the decomposition of $d\varphi$ we see that
\begin{align}\label{yayyo}
 (X\lrcorner d\varphi) \wedge \ast f\varphi &= \tau_0 f (X\lrcorner\ast\varphi)\wedge\ast\varphi + 3
 f \big(X\lrcorner(\tau_1\wedge\varphi)\big)\wedge\ast\varphi + f(X\lrcorner\ast\tau_3)\wedge\ast\varphi \notag\\
 &= 3
 f \big(X\lrcorner(\tau_1\wedge\varphi)\big)\wedge\ast\varphi + f(X\lrcorner\ast\tau_3)\wedge\ast\varphi\notag\\
 &= -3f(\tau_1\wedge\varphi)\wedge(X\lrcorner\ast\varphi) -
 f\ast\tau_3\wedge(X\lrcorner\ast\varphi)\notag\\
 &= -3f\tau_1\wedge\varphi\wedge(X\lrcorner\ast\varphi)\notag\\
 &= -3f\tau_1\wedge(-4\ast X^{\flat})\notag\\
 &= 12f\tau_1\wedge\ast X^{\flat},
\end{align}
where we have used the identity $\varphi\wedge(X\lrcorner\ast\varphi) = -4\ast X^{\flat}$ (see Appendix A in \cite{SK}) in 
the fifth equality and 
also the point-wise orthogonality of the $G_2$-decomposition of differential forms in the second and fourth equalities above.  
On the other hand, from the decomposition of $d\ast\varphi$ we have
\begin{align}\label{yay}
 \int_{M}d(X\lrcorner\varphi)\wedge\ast f\varphi &= \int_M
 (X\lrcorner\varphi)\wedge\ast\delta f\varphi\notag\\
 &= -\int_M (X\lrcorner\varphi)\wedge d\ast f\varphi\notag\\
 &= -\int_M (X\lrcorner\varphi)\wedge(df\wedge\ast\varphi +
 f\,d\ast\varphi)\notag\\
 &= -\int_M (X\lrcorner\varphi)\wedge df\wedge\ast\varphi
 - \int_M f(X\lrcorner\varphi)\wedge(4\tau_1\wedge\ast\varphi +
 \ast\tau_2)\notag\\
 &= -\int_M df\wedge\ast\varphi\wedge(X\lrcorner\varphi) -4\int_M f(X\lrcorner\varphi)\wedge \tau_1\wedge\ast\varphi
 \notag\\
 &= -\int_M df\wedge\ast\varphi\wedge(X\lrcorner\varphi) -4\int_M
 f\tau_1\wedge\ast\varphi\wedge(X\lrcorner\varphi)\notag\\
 &= -\int_M df\wedge 3\ast X^{\flat} - 4\int_M f\tau_1\wedge 3\ast X^{\flat}\notag\\
 &= -3\int_M df\wedge\ast X^{\flat}-12\int_M f\tau_1\wedge \ast X^{\flat},
\end{align}
where we have also used the identity
$\ast\varphi\wedge(X\lrcorner\varphi) = 3\ast
X^{\flat}$ (see Appendix A in \cite{SK}).  Integrating (\ref{yayyo}) and adding to (\ref{yay}), the lemma now follows. \qed \\

\newtheorem{corollary}{Corollary}
\begin{corollary}\label{main1}
 There are no compact expanding solitons,
 and there are no compact steady solitons except torsion-free $G_2$-structures.
\end{corollary}

\text{\bf{Proof.}} Wedging both sides of (\ref{sol}) by $\ast\varphi$ and integrating,
we have
\begin{align}\label{return}
  \rho\int_M \varphi\wedge\ast\varphi &= -\int_M
  \Delta\varphi\wedge\ast\varphi \notag\\
  &= -\int_M d\varphi\wedge\ast d\varphi \,-\, \int_M \delta\varphi \wedge\ast \delta\varphi  ,
\end{align}
where we have used Lemma \ref{wonder} with $f \equiv 1$.  Then since
$\int_M \varphi\wedge\ast\varphi$ is the volume and hence
non-zero, $\rho \leq 0$ necessarily because the right-hand side of (\ref{return}) is
non-positive.  The case of $\rho = 0$ is equivalent to $d\varphi = \delta\varphi = 0$ by (\ref{return}), hence torsion-free.
\qed  \\

As an offshoot to the proof of Corollary \ref{main1}, we also have the following observation.

\begin{corollary}
 For any compact soliton $\varphi$ satisfying (\ref{sol}), the constant $\rho$ 
 has the Rayleigh quotient expression: 
 \begin{equation}\label{quotient}
  \rho = \frac{-\int_M \Delta\varphi\wedge\ast\varphi}{\int_M \varphi\wedge\ast\varphi}.
 \end{equation}
\end{corollary} 

Equation (\ref{quotient}) means that for a compact soliton, $\rho$ is completely determined by the $G_2$-structure.  
This simple 
observation will play a role in Section \ref{rigid}.  Also recall from Section \ref{rambling} that $\rho = \dot{\tau}/\tau^{1/3}$ for any time $t$ within solution (\ref{soliton})'s existence.  Then (\ref{quotient}) implies that $\dot{\tau}$ will always have the same sign at any time $t$ within the soliton solution's existence, i.e. this means that 
a compact soliton solution is either "always" shrinking or "always" steady.

\section{Eigenforms as Shrinking Solitons}
 The defining equation (\ref{sol}) for a soliton shows that
it is a kind of generalized eigenvalue equation for $\varphi$
(with respect to $g_\varphi$).  In particular, when $\rho < 0$ an
\text{\it{eigenform}}: $-\Delta\varphi = \rho\varphi$ solves the
soliton equation with $X =0$ ($f_t = \text{Id} \,\,\forall t$).
Eigenforms as Laplacian solitons are analogous to Einstein
metrics as trivial examples of Ricci solitons.  Nevertheless, it
is enlightening to write out the exact solution for an eigenform
in terms of (\ref{soliton}) in the shrinking case.

\newtheorem{proposition}{Proposition}
\begin{proposition}
 Suppose $\varphi$ is an eigenform then it is a shrinking soliton
 with a solution in the form of (\ref{soliton}) as
 \begin{equation}\label{shrinker}
  \varphi(t) = \Big(1+\frac{2}{3}\rho t\Big)^{3/2}\varphi.
 \end{equation}
\end{proposition}

\text{\bf{Proof.}} We assume the solution is of the form
$\varphi(t) = R(t)\varphi$ for some real-valued function $R(t)$,
and we set $R(0) = 1$.  Then we see that
\[
  R'(t)\varphi = \frac{\partial\varphi(t)}{\partial t} =
  -\Delta_{\varphi(t)}\varphi(t) =
  -R(t)^{1/3}\Delta_{\varphi}\varphi = \rho R(t)^{1/3}\varphi.
\]
From this we must have $R'(t) = \rho R(t)^{1/3}$, which by
dividing both sides by $R(t)^{1/3}$ we can rewrite as
\begin{equation}\label{interlude}
 \frac{3}{2}\frac{d}{dt}\Big(R(t)^{2/3}\Big) = \rho.
\end{equation}
By integrating both sides of (\ref{interlude}) and using the
initial condition, we get exactly the
desired solution (\ref{shrinker}). \qed \\

From (\ref{shrinker}) we see that for an eigenform the
singularity time is $t = -3/2\rho >0$.  There are no eigen-forms
for the expanding case, thus if we want to find expanding solitons (of course,
only when $M$ is noncompcat) we must solve (\ref{sol}) with some
nontrivial vector field $X$.  A natural question is when is a closed $G_2$-structure an
eigenform. Noting the characterization of $\Omega^3_1$, we in
fact have the following result.

\begin{proposition}
 Let $\varphi$ be a closed $G_2$-structure.  Then $\varphi$ is an
 eigen-form if and only if $\Delta\varphi \in \Omega^3_1$.
\end{proposition}

\text{\bf{Proof.}} The only if part is trivial.  For the other
direction, note that if $\Delta\varphi = f\varphi$ then since
$d\varphi =0$ we must have $d(f\varphi) =0$ as well.  In other
words,
\[
  df\wedge\varphi = 0.
\]
Recall that the equation above, along with the special form that
$\varphi$ takes with respect to a local orthonormal frame, shows
in a straight-forward way that $f$ must be constant.  Thus
$\varphi$ must be an eigenform.  \qed \\

Although we have the proposition above, it is not
straight-forward to find closed eigenforms.  However, there are
plenty of eigenforms that are not closed.  We recall the following

\begin{definition}
 A $G_2$-structure is called nearly parallel if its only
 nonzero torsion form is $\tau_0$.  We say that a $7$-manifold is
 a nearly $G_2$ manifold if it admits a nearly parallel $G_2$-structure.
\end{definition}

Thus $(M,\varphi)$ is a nearly $G_2$ manifold if and only if
$\delta\varphi = 0$ and $d\varphi = \tau_0\ast\varphi$.  We want
to point out that in this case, $\tau_0$ is necessarily a constant
because $0 = d(d\varphi) = d\tau_0\wedge\ast\varphi$ implies so.
The squashed $7$-sphere is an example of a nearly $G_2$ manifold.
By straight-forward computation, we immediately see that:

\begin{proposition}
 If $(M,\varphi)$ is a nearly $G_2$ manifold, then $\varphi$ is an
 eigenform satisfying
\begin{equation}
 \Delta\varphi = \tau_0^2\varphi.
\end{equation}
In particular, $\varphi$ is a shrinking soliton.
\end{proposition}

Thus we see that any nearly $G_2$ manifold admits a shrinking
soliton.\\

Let us return to the somewhat opposite case to a nearly $G_2$
manifold, which is the case where a $G_2$-structure $\varphi$ is
closed.  In this case, we can consider the de Rham cohomology
classes.  We would first like to point out that if $\varphi$ is
closed then $f_t^*\varphi$ always stays within the cohomology
class of $f_0^*\varphi$ ($\varphi$ is closed, hence 
$f_0^*\varphi$ is too),\footnote{Clearly, a soliton solution remains
closed for all $t$ if $\varphi$ is closed.  On the other hand, in
general one can show from (\ref{lap}) that it preserves the
closedness condition of the initial value $\varphi(0)$.} which is 
based on the following elementary result:

\begin{lemma}
Let $\omega$ be a closed $p$-form, then for any smooth family of
diffeomorphisms $f_t$ homotopic to the identity, $f_t^*\omega$ and
$\omega$ are cohomologous for all $t$.
\end{lemma}

\noindent Now, for a smooth family $f_t$ of
diffeomorphisms of $M$, we define a new family $\tilde{f}_t =
f_0^{-1}\circ f_t$.  Then $\tilde{f}_t^*\circ f_0^*\varphi =
f_t^*\varphi$, and $\tilde{f}_0 = I$.  Thus by the lemma above, we
see that $f_t^*\varphi$ is cohomologous to $f_0^*\varphi$ for all
$t$.  Therefore a closed\footnote{Note that we are distinguishing
this from compact solitions.} steady soliton solution remains in the
original cohomology class of $f_0^*\varphi$ (normalizing the scaling factor to be 
$1$) and hence can be
seen as a periodic solution in a fixed cohomology class.  On the other hand,
for closed expanding and shrinking solitons, we see that
\begin{equation}
 \rho\varphi = -d(X\lrcorner\varphi + \delta \varphi).
\end{equation}
In other words, the $G_2$-structure $\varphi$ must be
exact.\footnote{If $M$ is compact, this immediately precludes a
closed shrinking
soliton from being torsion-free.  See the next section for a more general derivation of this fact.}\\

\section{Rigidity of Laplacian Solitons}\label{rigid}
On a compact $7$-manifold, if a $G_2$-structure $\varphi$ satisfies the soliton equation
(\ref{sol}) for some $X$ and $\rho$, we may ask whether there
are other vector fields $X'$ and constants $\rho'$ with which
$\varphi$ is also a soliton.  We already saw at the end of Section \ref{cpt} that 
we must have $\rho' = \rho$.  On the other hand, if $-\Delta\varphi = L_{X'} \varphi + \rho\varphi$ for 
some other vector field $X'$, then subtracting it from the original soliton equation gives 
\[
 L_{X-X'}\varphi = 0.
\]
In other words, for any compact soliton $\varphi$, the
only symmetries of its defining equation (\ref{sol}) are $X
\longrightarrow X + Y$ for vector fields such that $L_Y\varphi
=0$. With the $G_2$-structure fixed, this is the only change to a soliton equation that 
leaves it invariant.   In general, a vector field $X$ such that
$L_X\varphi = 0$ is simply called a \text{\it{symmetry}} of the $G_2$-structure
$\varphi$.  One can go further in revealing the properties of these
vector fields, in fact we have the following result.

\begin{proposition} \label{divfree}
 On a compact $7$-manifold $M$ admitting a $G_2$-structure $\varphi$, any symmetry $X$ of $\varphi$ must satisfy
 $\text{div}(X) =0$.
\end{proposition}

\text{\bf{Proof.}} Again by Lemma \ref{wonder}, if $L_X\varphi =
0$, then
\begin{align}
  0 = \int_M L_X\varphi\wedge\ast f\varphi &= -3\int_M df\wedge\ast
   X^{\flat} \notag\\
   &= -3\int_M f\,\ast\delta X^{\flat}
\end{align}
for all $f\in C^{\infty}(M)$.  This implies that $\delta X^{\flat}
= 0$, or $\text{div}(X) =0$. \qed \\

  Note that Proposition \ref{divfree} is a general
result for any $G_2$-structure $\varphi$ on a compact $7$-manifold. \\

We want to understand the full structure of the space of symmetries for solitons in general.
However, for now we will prove a partial result, but which has
immediate significance.

\begin{lemma}\label{1dir}
   If a $G_2$-structure $\varphi$ is closed, then
 $X\lrcorner\varphi$ is harmonic for any symmetry $X$ of
 $\varphi$.  Moreover, if $M$ is compact and $\varphi$ is torsion-free then $X^{\flat}$
 is a harmonic $1$-form, and in particular if $g_\varphi$ has
 full $G_2$ holonomy then $X=0$.
\end{lemma}

\text{\bf{Proof.}} If $d\varphi = 0$, then we have
\[
   L_X\varphi = d(X\lrcorner\varphi) = 0.
\]
Then using the fact that $X\lrcorner\varphi \in \Omega^2_{7}$,
\begin{align}\label{2birds}
   d\ast(X\lrcorner\varphi) &=
   -\frac{1}{2}d\big[\varphi\wedge(X\lrcorner\varphi)\big)\notag\\
   &= -\frac{1}{2}\big(d\varphi\wedge(X\lrcorner\varphi) - \varphi\wedge
   d(X\lrcorner\varphi)\big]\notag\\
   &= 0.
\end{align}
Thus $X\lrcorner\varphi \in \Omega^2_{7}$ is
harmonic.\footnote{Note that the computations in (\ref{2birds})
show that if $\varphi$ is closed, then $(X\lrcorner\varphi)$ being closed
implies $(X\lrcorner\varphi)$ is coclosed as well, for any vector field $X$.}\\

Now, if $\varphi$ is torsion-free then we further have
\begin{align}
 0 &= d\ast(X\lrcorner\varphi) \notag\\
 &= d(X^{\flat}\wedge\ast\varphi) \notag\\
 &= dX^{\flat}\wedge\ast\varphi. \notag
\end{align}
 This shows that $d X^{\flat} \in \Omega^2_{14}$.  Then using the other characterization of $\Omega^2_{14}$ we see that 
\begin{align}
 d\big(X^{\flat}\wedge\varphi\wedge d X^{\flat}\big) &= dX^{\flat}\wedge\varphi\wedge dX^{\flat}\notag\\
 &= dX^{\flat}\wedge\ast dX^{\flat} \notag
\end{align}
Integrating both sides above and using Stoke's theorem, we see that $dX^{\flat} =0$ necessarily.  
By Proposition \ref{divfree}, we conclude
that $X^{\flat}$ is harmonic.  If in addition we have full $G_2$
holonomy, then we know that $H^1(M,\mathbb{R}) = 0$\footnote{For example, see \cite{DJ}.} and so
$X^{\flat} =0$. \qed \\

\begin{corollary}\label{isom}
 Let $\varphi$ be a torsion-free $G_2$-structure on a compact
 $7$-manifold $M$.  Then a vector
 field $X$ is a symmetry of $\varphi$ if and only if $X^{\flat}$
 is a harmonic $1$-form.  Thus the space of symmetries of
 $\varphi$ is isomorphic to $H^1(M,\mathbb{R})$.
\end{corollary}

\text{\bf{Proof.}}  The only if part is given by the lemma above.
To prove the if part, we will employ the general relation below
for the Levi-Civita connection:
\begin{align}\label{realized}
 L_X\varphi (Y_1,Y_2,Y_3) &= \nabla_X\varphi\,(Y_1,Y_2,Y_3) +
 \varphi(\nabla_{Y_1}X, Y_2,Y_3) +
 \varphi(Y_1,\nabla_{Y_2}X,Y_3) \notag\\
 &\,\,\,\,+
 \varphi(Y_1,Y_2,\nabla_{Y_3}X)
\end{align}
for any vector fields $X,Y_1,Y_2,Y_3$.  The torsion-free condition is defined by
$\nabla\varphi = 0$.  Furthermore, since any $G_2$ manifold must have zero Ricci
curvature everywhere, by the Bochner's Theorem we know that any harmonic
$1$-form must be parallel.  Then we must have $\nabla X = 0$.  Combining these facts
into (\ref{realized}) we get the desired result.  \qed \\

Corollary \ref{isom} contains the special case stated at the end
of Lemma \ref{1dir}.  In general, we know that for a compact $G_2$
manifold the following holds:
\begin{align}\label{diag}
 Hol(M) = \{1\} &\Longleftrightarrow b_1(M) = 7\notag\\
 Hol(M) = SU(2) &\Longleftrightarrow b_1(M) = 3\notag\\
 Hol(M) = SU(3) &\Longleftrightarrow b_1(M) = 1\notag\\
 Hol(M) = G_2 &\Longleftrightarrow b_1(M) = 0.
\end{align}
Thus Corollary \ref{isom} shows that for a torsion-free $G_2$
structure, increasing the holonomy will decrease its symmetries -
making it more and more "rigid".  In particular, the $G_2$
structure of a $G_2$ manifold with full $G_2$ holonomy admits no
non-trivial symmetries, hence it is "rigid".  This is a rigidity
result in stark contrast to that for compact manifolds with
positive Ricci curvature, whose $b_1(M) = 0$ but yet
typically admits a lot of Killing vector fields. \\

For $G_2$ manifolds that do not have full $G_2$ holonomy, the list
(\ref{diag}) shows that there are non-trivial symmetries $X$.  In
other words, such $X$ generates a one-parameter family $f_t$ of
diffeomorphisms homotopic to the identity such that $f_t^*\varphi
= \varphi$ for all $t$.  The existence of such symmetries is most visible in
the following known examples:

\begin{enumerate}
 \item $\mathbb{T}\times Y$, with $\varphi = dx\wedge\omega + Re\,\theta$, where
 $Y$ is a Calabi-Yau $3$-fold with K\"{a}hler form $\omega$ and
 holomorphic volume form $\theta$.  The holonomy group is $SU(3)$.

 \item $\mathbb{T}^3\times Y$, with $\varphi = dx_{123} + dx_1\wedge\omega +
 dx_2\wedge \text{Re}\,\theta -dx_3\wedge \text{Im}\,\theta$, where $Y$ is a
 Calabi-Yau $2$-fold with K\"{a}hler form $\omega$ and
 holomorphic volume form $\theta$.  The holonomy group is $SU(2)$.

 \item $\mathbb{T}^7$, inheriting the standard $G_2$ structure on
 $\mathbb{R}^7$.  The holonomy group is $\{1\}$.
\end{enumerate}

We would like to point-out that Corollary \ref{main1} and Corollary \ref{isom} together say that a compact steady soliton $(\varphi, X)$ consists of a 
torsion-free $G_2$-structure $\varphi$ and a symmetry of $\varphi$.     
As mentioned in the Introduction, it would be desirable to prove an analogous result for 
compact shrinking solitons.  However, it seems that there are some technical difficulties.

\section{Concluding Remarks and Questions}
In this paper we have discussed fundamental properties of Laplacian solitons on manifolds admitting a $G_2$-structure.  
In particular, we have investigated solitons on compact $7$-manifolds to some detail.  Due to the recent flurry of 
interest in Ricci solitons, we want to make some comparisons to it.  Recall that a Ricci soliton is a metric 
$g$ along with a vector field $X$ on a manifold $M$ such that 
\begin{equation}\label{rs}
 \text{Ric}_g = L_X g + \rho g
\end{equation}
for some constant $\rho$.  The sign on $\rho$ dictates the terminologies of shrinking, steady, and expanding solitons 
in the same way as for our Laplacian solitons.  \\

When $M$ is compact, we know that the only steady and expanding 
Ricci solitons are Einstein, i.e. $X$ must be a Killing vector field.  Corollary \ref{main1} can be seen as a 
close parallel to this result.  In fact, by the works of Richard Hamilton and Thomas Ivey we also know that 
compact shrinking Ricci solitons in dimensioins $2$ and $3$ have to also be Einstein.  Since we mentioned that 
eigenforms should be viewed as an analogy to Einstein metrics in the context of solitons, it is a natural 
question whether or not compact shrinking Laplacian solitons also have to be eigenforms.  If this conjecture 
is true, then we would be able to finish the classification of compact Laplacian solitons if we further 
classify all compact eigenforms.\footnote{A relevant question here would be: are the only compact eigenforms 
nearly parallel $G_2$-structures?}  The present 
difficulty seems to be due to a lack of maximum principle-type techniques in $G_2$ geometry, since such techniques 
were essential in proving results in Ricci Solitons (or the Ricci flow in general).      \\

Using his entropy functionals, when $M$ 
is compact G. Perelman showed that in (\ref{rs}) $X$ can always be chosen as a gradient vector field, i.e. compact Ricci solitons 
are always \text{\it{gradient}} solitons.  It would be interesting to know if compact Laplacian solitons are always 
gradient as well, in the same sense.  The entropy functionals of Perelman has a deeper implication: they turn the Ricci 
flow into a gradient-like\footnote{coupled to other equations} flow with respect to the functionals.  In particular, compact Ricci solitons appear as 
critical points of these functionals.  Functionals 
associated to the Laplacian flow have been suggested in \cite{NH} and \cite{WW}, and it would be nice to see if they 
(or possibly other functionals) can be used to characterize compact Laplacian solitons.  \\

Finally, recall that on a compact manifold $M$ of nonpositive Ricci curvature the space of Killing vector fields is contained 
in the space of parallel vector fields, and if the Ricci curvature is negative there are no nontrivial Killing vector fields.\footnote{Therefore when $M$ is Ricci-flat the space of Killing vector fields is isomorphic 
to a subspace of $H^{1}(M,\mathbb{R})$ via Bochner's Theorem.}  In view of the preceding results and our intended analogy, 
it seems that for compact shrinking Laplacian solitons all symmetries should be parallel vector fields (with respect 
to $g_\varphi$) as well, or even always trivial.  On the other hand, it seems more tantalizing to conjecture that the space of symmetries of a 
compact shrinking Laplacian soliton is simply the eigenspace of the equation $\Delta \omega = \rho\omega$ for $1$-forms 
$\omega$.  Identifying the symmetries of solitons will be important because one would like to 
construct the associated \text{\it{moduli spaces}} by dividing out the symmetries.\footnote{A discussion on such moduli spaces 
have also appeared in \cite{WW2}.}

\vskip 1cm

\small
\textsc{Department of Mathematics, Case Western Reserve University, Cleveland OH   44106}\\
E-mail address, Christopher Lin: \texttt{ccl37@case.edu}


\begin{thebibliography}{10}

\bibitem{RB} Robert Bryant, \emph{Some remarks on $G_2$-structures}, Proceedings of Gokova Geometry/Topology Conference, 
 Gokova, 2006, pp. 75-109.
\bibitem{BX} Robert Bryant and Feng Xu, \emph{Laplacian flow for closed $G_2$-structures: Short time behavior}, arXiv:1101.2004
\bibitem{FG} M. Fern\'{a}ndez and A. Gray, \emph{Riemannian manifolds with structural group $G_2$}, Ann. Mat. Pura Appl 
(IV) \text{\bf{32}} (1982), 19-45. 

\bibitem{NH} Nigel Hitchin, \emph{The geometry of three-forms in six and seven dimensions}, arXiv:math/0010054

\bibitem{DJ} Dominic D. Joyce, \emph{Compact manifolds with special holonomy}, Oxford Mathematical Monographs, 
Oxford University Press, Oxford, 2000. MR MR1787733 (2001k:53093)



\bibitem{SK} Spiro Karigiannis, \emph{Flows of $G_2$-structures I}, Q.J. Math. \text{\bf{60}} (2009), no.4, 487-522.

\bibitem{KMT} Spiro Karigiannis, Benjamin McKay, Mao-Pei Tsui, \emph{Soliton solutions for the 
Laplacian coflow of some $G_2$ structures with symmetry}, arXiv:1108.2192v1
\bibitem{WW} Hartmut Weiss and Frederik Witt, \emph{A heat flow for special metrics}, arXiv:0912.0421 

\bibitem{WW2} Hartmut Weiss and Frederik Witt, \emph{Energy Functionals and Soliton Equations for $G_2$ forms}, 
arXiv:1201.1208v1

\bibitem{XY} Feng Xu and Rugang Ye, \emph{Existence, convergence and limit map of the Laplacian Flow}, arXiv:0912.0074
\end{thebibliography}
\end{document}